\documentclass{amsart}

\usepackage{amssymb}
\usepackage{verbatim}

\newtheorem{Th}{Theorem}[section] 
\newtheorem{Lem}[Th]{Lemma} \newtheorem{Prop}[Th]{Proposition}

\newtheorem{Rem}[Th]{Remark}

\newtheorem{claim-num}{Claim}

\newtheorem*{theo}{Theorem}
\newtheorem*{cor}{Corollary}

\numberwithin{equation}{section}
\renewcommand{\theequation}{\thesection.\arabic{equation}}

\def\aut#1{\operatorname{Aut}(#1)}

\def\gl#1{\operatorname{GL}(#1)}
\def\spl#1{\operatorname{SL}(#1)}
\def\End#1{\operatorname{End}(#1)}

\def\id{\operatorname{id}}
\def\inv{^{-1}}
\def\rank{\operatorname{rank}}
\def\str#1{\langle#1\rangle}

\def\f{\varphi}

\def\s{\sigma}
\def\vk{\varkappa}
\def\Theo{\operatorname{Th}}

\def\av#1{\overline{#1}}

\def\cB{{\mathcal B}}

\def\cD{{\mathcal D}}

\def\N{\mathbf N}
\def\Z{\mathbf Z}

\renewcommand{\le}{\leqslant}
\renewcommand{\ge}{\geqslant}

\begin{document}


\title[Automorphism groups of free abelian groups]
{What does the automorphism
group \\ of a free abelian group $A$ know
about $A$?}
\author{Vladimir Tolstykh}
\thanks{Supported by a NATO PC-B grant via The
Scientific and Technical Research Council of Turkey (T\"UB\'ITAK)}
\address{Department of Mathematics\\ Yeditepe University\\
34755 Kay\i\c sda\u g\i \\
Istanbul\\
Turkey}
\maketitle

\begin{abstract}
Let $A$ be an infinitely generated free abelian group.
We prove that the automorphism group $\aut A$
first-order interprets the full second-order theory of
the set $|A|$ with no structure. In particular, this
implies that the automorphism groups of two infinitely
generated free abelian groups $A_1,A_2$ are
elementarily equivalent if and only if the sets
$|A_1|,|A_2|$ are second-order equivalent.
\end{abstract}

\section*{Introduction}

In his paper \cite{ShCler} of 1976 Shelah proved that
the elementary theories of the endomorphism
semi-groups of free algebras of `large' infinite ranks
had very strong expressive power. More precisely, let
$\bold V$ be an {\it arbitrary} variety of algebras
and $F_\vk(\bold V)$ be a free algebra from $\bold V$ with $\vk
\ge \aleph_0$ free generators. Then the
endomorphism semi-group $\End{F_\vk(\bold V)}$ first-order
interprets the full-second theory $\Theo_2(\vk)$ of
the cardinal $\vk$ (viewed as a set with no
structure), provided that $\vk$ is greater
than the cardinality of the language of $\bold V.$

That remarkable result naturally leads to the
following problem: what are the varieties of algebras
for which the {\it automorphism groups} of free
algebras are logically strong in a similar sense? Shelah himself
formulated this problem in the cited paper \cite{ShCler}
and then after more than 20 years mentioned it again in
his survey \cite{ShMatJap}: Problem
3.14 from \cite{ShMatJap} suggested to classify the varieties of algebras
$\bold V$ such that the automorphism groups
$\aut{F_\vk(\bold V)}$ first-order interpret
the theory $\Theo_2(\vk)$ for all
(or all sufficiently large) infinite cardinals $\vk.$

 The results
on symmetric groups obtained by Shelah before the
publication of the paper \cite{ShCler} implied that,
for instance, the variety of all sets with no
structure and the variety of all semi-groups were the
examples of, say, `negative' kind. Indeed, according
to \cite{Sh1}, the symmetric group of an infinite
cardinal $\vk,$ in other words, the automorphism
group of the set $\vk$ with no structure,
first-order interprets the theory $\Theo_2(\vk)$
only if the cardinal $\vk$ is `small' (namely,
at most $2^{\aleph_0}$).

The author found in \cite{ToAPAL}--as a byproduct of
his study of the elementary types of infinite-dimensional
classical groups--that for any variety of vector spaces
the automorphism groups of free algebras are as logically
strong as the endomorphism semi-groups. A bit
informally, one of the results from \cite{ToAPAL} can be
quoted in the following form: if $\vk$ is an infinite
cardinal, then the general linear group $\gl{\vk,D}$
over a division ring $D$ first-order interprets
$\Theo_2(\vk),$ provided that $\vk > |D|.$ Thus
varieties of vector spaces give examples of `positive'
kind as to Shelah's problem.

In the papers \cite{ToJLM2} and \cite{ToContMat} the
author studied Shelah's problem for classical group
varieties.  It turned out that the variety of all
groups and any variety $\frak N_c$ of nilpotent groups
of class $c \ge 2$ meet the requirements of Shelah's
problem: if $F$ is an infinitely generated free or
free nilpotent group, then the group $\aut F$
first-order interprets the theory $\Theo_2(|F|)$ $(=\Theo_2(\rank F).)$  In
the present paper we examine the case of the variety
of all abelian groups.  The main result of the paper
states that the variety in question also meets
requirements of Shelah's problem.

Let $A$ denote an infinitely generated free
abelian group; clearly, $A$ can be considered as a
free $\Z$-module. One of the standard approaches to understanding of the
nature of the automorphism groups of modules is an
investigation of possibility of generalization for
these groups of the methods developed for general linear
 groups, the automorphism groups of vector
spaces.   In the first section of the paper
we, like in \cite{ToAPAL}, work to reconstruct by
means of first-order logic in $\aut A$ some geometry
of the $\Z$-module $A.$ Namely, we interpret in $\aut A$
the family $\cD^1(A)$ consisting of all direct summands
of $A$ having rank or corank one. To make comparison, the
first-order interpretation in the general linear group
$\gl V$ of an infinite-dimensional vector space $V$
of the family of all lines and hyperplanes of $V$ done in
\cite{ToAPAL} is much longer. However, both
interpretations have much in common and
both originated from the well-known
works on classical groups.

In principle, the reconstruction of $\cD^1(A)$ can be
extended to the reconstruction in $\aut A$ of the family
$\cD(A)$ of all direct summands of $A$ followed by the
first-order interpretation in the structure $\str{\aut A,\cD(A)}$ of the
endomorphism semi-group $\End A$ of $A$ (similarly to
\cite{ToAPAL}). We, however, prefer a shorter way,
making in Section 2 an effort to reconstruct in $\aut
A$ the general linear group of some vector space
of dimension $|A|.$ Namely, using the action of
$\aut A$ on $\cD^1(A)$ we prove
$\varnothing$-definability in $\aut A$ of the
principal congruence subgroup $\Gamma_2(A)$ of level
two. The quotient subgroup $\aut A/\Gamma_2(A)$ is
isomorphic to the general linear group of the vector
space $A/2A$ over the field $\Z_2.$ Thus
the group $\aut A$ first-order interprets
the group $\gl{|A|,\Z_2}.$ The latter group,
as it has been said above, first-order interprets
the theory $\Theo_2(|A|).$
As a consequence, we have that
the automorphism groups $\aut{A_1}$
and $\aut{A_2},$ where $A_1,A_2$
are infinitely generated free
abelian groups, are elementarily
equivalent if and only the cardinals
$|A_1|$ and $|A_2|$ are second-order
equivalent as sets.

The author is very grateful to Oleg Belegradek for his
kind and genuine interest in this research and for his
valuable comments on the first draft of this paper.
The author would like also to thank Valery Bardakov
for helpful discussions.

\section{Definable geometric properties of automorphisms}

Let $A$ denote a free abelian group of infinite rank.
As it has been said in the Introduction, our aim in
this section is a first-order reconstruction in $\aut
A$ of the family of direct summands of $A$ of rank or
corank one (we say that a direct summand $B$ of $A$
has {\it corank $m$}, if any direct complement of $B$
to $A$ is of rank $m$.)

We shall essentially exploit the structure of
involutions (the elements of the order two) in the
group $\aut A$ given by the following
theorem.

\begin{Th} \label{HRCanForm}
Let $G$ be a free abelian group. Every
involution $\f \in \aut G$ has a basis
$\cB$ of $G$ such that for any $b \in \cB$
either $\f b=\pm b,$ or $\f b \in \cB.$
\end{Th}

The theorem was first established for the groups of
finite rank by Hua and Reiner \cite[Lemma 1]{HuaRei};
in general, the result is proven in \cite{ToCamb}.
Let us call a basis of $A$ on which $\f$ acts in a way
described in the Theorem a {\it canonical} basis
for $\f.$

Let $2A$ denote the group of even
elements of $A$:
$$
2A=\{2a : a \in A\}.
$$
The natural homomorphism $A \to A/2A$
induces the homomorphism of the
automorphism groups $\aut A \to \aut{A/2A}$
which we will denote by $\widehat{\phantom a}.$
The fact that the group $A/2A$ can be viewed
as a vector space over $\Z_2$ will be extensively
used in this paper.

\begin{Rem} \label{To-Can-Forms}
\em Take an involution $\f \in \aut A$ and some its
canonical basis $\cB.$ The (cardinal) number $p(\cB)$
of unordered pairs $\{b,\f b\},$ where $b \in \cB$
and $\f b \ne \pm b$ is an invariant of $\f.$ Indeed,
$p(\cB)$ equals the {\it residue} of the
induced linear transformation $\widehat{\f}$
of the vector space $A/2A$ over $\Z_2$:
$$
p(\cB) =\text{res}(\widehat \f) =\dim \text{Res}(\widehat \f).
$$
(here $\text{Res}(\widehat \f)$ is the
image of the linear transformation $1-\widehat \f,$
see \cite{O'M}). This implies that if $(\f_1,\cB_1),$
$(\f_2,\cB_2)$ are pairs similar to the pair
$(\f,\cB)$ and $\f_1,\f_2$ are conjugate in $\aut A,$
then $p(\cB_1)=p(\cB_2).$ \end{Rem}

Let $\f$ be an involution in $\aut A$; we let $A^+_\f$ and $A^-_\f$
denote the subgroups
$$
\{a : \f a =a\} \text{ and } \{a: \f a = -a\}
$$
respectively; clearly, $\f$ is diagonalizable
if and only
$$
A = A^+_\f \oplus A^-_\f.
$$
It is helpful to remember that two diagonalizable
involutions from $\aut A$ are commuting if and only
if there is a basis of $A$ in which they both
diagonalizable.

We shall call a diagonalizable involuton $\f$ a $\gamma$-{\it involution}, where $\gamma$ is
a cardinal, if
$$
\gamma = \rank A^-_\f < \rank A^+_\f.
$$
$1$-involutions, like in linear group theory, will be
called {\it extremal} involutions.

A number of facts on definability of certain families
of involutions in the automorphism groups of
infinitely generated free abelian groups has been
proved implicitly in the author's paper \cite{ToCamb}.
Because of that we shall give only sketches of proofs
for the next two statements, Lemma \ref{Diags-Are-Def}
and Lemma \ref{1-2-4-Def}; the reader is
referred to the proof of Proposition 2.4  in
\cite{ToCamb} to find there the omitted details.

For an involution $\f$ in the group $\aut A$ we shall
denote by $K(\f)$ the conjugacy class of $\f$ in $\aut A.$
The set $K^2(\f)=K(\f)K(\f)$ is the family of all
products $\f_1\f_2,$ where $\f_1,\f_2 \in K(\f).$

\begin{Lem} \label{Diags-Are-Def}
The family of all diagonalizable involutions
is $\varnothing$-definable in \linebreak $\aut A.$
\end{Lem}

\begin{proof}
We claim that $\f$ is diagonalizable if and only
if the set $K^2(\f)$ contains no elements of
order three.

Using Theorem \ref{HRCanForm} one checks that the
diagonalizable involutions are exactly involutions in
the kernel of the homomorphism $\widehat{\phantom a} : \aut A \to
\aut{A/2A}.$ On the other hand, the images under $\widehat{\phantom a}$
of all elements of order three from $\aut A$ are
non-trivial. This implies that if $\f$ is
diagonalizable, then there are no elements of order
three in $K^2(\f).$

Conversely, for any non-diagonalizable involution
$\psi \in \aut A$ we can easily find a conjugate $\psi'$
of $\psi$ such that the automorphism $\psi\psi'$ is of
order three.  \end{proof}

\begin{Lem} \label{1-2-4-Def}
The families of extremal involutions
{\em(}$1$-involutions{\em)}, $2$-involutions
and $4$-involutions are all $\varnothing$-definable in $\aut A.$
\end{Lem}

\begin{proof}
A diagonalizable involution $\f$ is an extremal
involution if and only if all involutions
in $K^2(\f)$ are conjugate and $\f$ is not a
square in $\aut A.$

Indeed, if $\f$ is an extremal involution, then the
only involutions in the set $K^2(\f)$ are
$2$-involutions.  In particular, all involutions in
$K^2(\f)$ are conjugate. Applying Theorem
\ref{HRCanForm}, we can demonstrate that the latter
property holds also only for diagonalizable
involutions $\rho$ such that
$$
\rank A^+_\rho =1.
$$
But any such an involution is a square
in $\aut A,$ whereas any $1$-involution
is not.

The $2$-involutions are the only involutions from
$K^2(\f),$ where $\f$ is an arbitrary $1$-involution.
Let $\theta$ be a $2$-involution. Then $4$-involutions
are those involutions in $K^2(\theta)$ that are not
conjugate to $\theta.$
\end{proof}

We need also a family of {\it non-diagonalizable}
involutions $\{\pi\}$ whose elements satisfy
the condition
\begin{equation}
\rank A^+_\pi=1 \text{ or } \rank A^-_\pi=1.
\end{equation}
For any canonical basis $\cB$ for a non-diagonalizable
involution $\pi$ with (\theequation) we have that
\begin{itemize}
\item[(a)] $\cB$ contains exactly one pair of distinct elements,
say, $b,c$ taken by $\pi$ to one another (Remark
\ref{To-Can-Forms});
\item[(b)] $\pi$ either inverts {\it all} elements
in $\cB\setminus \{b,c\},$ or fixes
{\it all} these elements (otherwise, both
subgroups $A^+_\pi$ and $A^+_\pi$
were of rank $> 1$).
\end{itemize}

Thus either $\pi \sim \pi',$ or $\pi \sim -\pi'$
for every pair of non-diagonalizable involutions $\pi,\pi'$
with (\theequation), where $\sim$ denotes
the conjugacy relation. Keeping in mind (a), we shall
call non-diagonalizable involutions with
(\theequation) by $1$-{\it permutations.}

\begin{Lem}
The following statements are equivalent:
\begin{itemize}
\item[(i)] $\pi$ is a $1$-permutation;
\item[(ii)] $\pi$ is not diagonalizable and the set $K^2(\pi)$ contains
no $4$-involutions.
\end{itemize}
In particular, the family of $1$-permutations
is $\varnothing$-definable in $\aut A.$
\end{Lem}

\begin{proof}
Let $\pi$ be a non-diagonalizable involution, which
is not a $1$-permutation, and let $\cB$ be a canonical basis for
$\pi.$ One then can readily find $\pi',$ a conjugate of
$\pi,$ whose product with $\pi$ is a $4$-involution.
Indeed, suppose first that $p(\cB) > 1$
(the notation was introduced in Remark
\ref{To-Can-Forms}). Then $\cB$ contains distinct
elements $b_1,b_2,b_3,b_4$ such that $\pi b_1=b_2$
and $\pi b_3=b_4.$

The second case is the case when $p(\cB)=1.$ Here
$\pi b_1=b_2$ for some distinct $b_1,b_2 \in \cB$
and, since $\pi$ is not a $1$-permutation, two such
elements $b_3$ and $b_4$ can be found in $\cB$ that
$$
\pi b_3 = b_3 \text{ and } \pi b_4 = -b_4.
$$

Then, for both of the cases under consideration, we
construct $\pi'$ as follows:  $\pi'b_i=-\pi b_i$ for
$i=1,\ldots,4$ and
$\pi' b=\pi b$ for all $b \in \cB\setminus
\{b_1,b_2,b_3,b_4\}.$

Conversely, suppose $\pi$ is a 1-permutation. We may
assume that $\rank A_\pi^-=1.$ Let $\pi_1, \pi_2$ be
conjugates of $\pi.$ Then $\operatorname{Im}(1-\pi_1)$
and $\operatorname{Im}(1-\pi_2)$ are subgroups of rank
1.  Since
$$
1-\pi_1\pi_2=(1-\pi_2)+(1-\pi_1)\pi_2,
$$
we have
$$
\operatorname{Im}(1-\pi_1\pi_2)\subseteq
\operatorname{Im}(1-\pi_1)+\operatorname{Im}(1-\pi_2),
$$
and so $\rank \operatorname{Im}(1-\pi_1\pi_2)\le 2.$
Then $\pi_1\pi_2$ is not a 4-involution because
for any 4-involution $\psi$ we have
$\rank \operatorname{Im}(1-\psi)=4.$
\end{proof}

Until the end of this section we fix some $2$-involution
$\theta^*.$ In order to mark somehow one special type of
commutativity with $\theta^*,$ we say that an extremal
involution $\psi$ (resp. a $1$-permutation $\psi$) commutes with
$\theta^*$ {\it properly}, if $\psi \sim \theta^*
\psi.$

We fix also an extremal involution $\f^*$ and
a $1$-permutation $\pi^*$ both properly
commuting with $\theta^*$ such that
$$
(\pi^* \f^*)^2 = \theta^*.
$$
Let $B$ denote the subgroup $A^-_{\theta^*}.$
Since both $\f^*$ and $\pi^*$ commute
with $\theta^*,$ they both preserve
$B$:
$$
\f^* B = \pi^* B=B.
$$
Since, further, $\f^*$ and $\pi^*$ commute
with $\theta^*$ properly, their
restrictions to $B$ are
an extremal involution and
a $1$-permutation of $\aut B,$
respectively.
Let
$$
f^* = \f^*|_B \text{ and } p^*=\pi^*|_B.
$$
We have that
$$
f^* p^* f^* p^* =-\id_B
$$
and then $p^* f^* p^*=-f^*.$ This implies
that $p^*$ takes to each other
the subgroups $A^+_{f^*}$ and $A^-_{f^*}$:
$$
p^* A^+_{f^*} = A^-_{f^*}.
$$
If then $e_1$ is a basis element
of $A^+_{f^*},$ then $e_2=p^* e_1$ is a basis
element of $A^-_{f^*}.$ Summing up, we
see that in the basis $\{e_1,e_2\}$ of $B$
the automorphisms $f^*$ and $p^*$ have
the matrices
$$
\begin{pmatrix}
1 & 0 \\
0 & -1
\end{pmatrix}
\text{ and }
\begin{pmatrix}
0 & 1\\
1 & 0
\end{pmatrix},
$$
respectively.

Next is the proof of $\varnothing$-definability of
certain transvections. Recall that a {\it unimodular}
element of $A$ (a {\it primitive} element in a more general
context) is an element of $A$ that can be included in
some basis of $A.$ Let $\delta : A \to \Z$ be a
non-zero homomorphism of abelian groups;
in this case the $\ker\delta$ is a direct summand of $A$ of corank 1.
Fix a unimodular element $x$ in $\ker \delta.$ Then the mapping
$$
\tau a =a +\delta(a) x,
$$
is an automorphism
of $A,$ which is called
a {\it transvection.} If $\tau$ is a transvection
determined by a homomorphism $\delta,$ then
one may correctly associate with $\tau$
a natural number defining it via
$$
m(\tau) = |\delta(y)|
$$
where $y \in A$ satisfies $A =\str y \oplus \ker
\delta.$ It can be easily seen that
for every pair $\tau_1,\tau_2$
of transvections $m(\tau_1)=m(\tau_2)$
if and only if $\tau_1$ and
$\tau_2$ are conjugate. We shall call a transvection
$\tau$ an $m$-{\it transvection}, if $m(\tau)=m.$

\begin{Lem} \label{4-1-Perms}
{\em (i)} Among the conjugates
$\rho$ of $\pi^*$ properly
commuting with $\theta^*$ there are
exactly four ones different
from $\pi^*$ that satisfy the equation
$$
(\pi^* \rho)^3 =\id_A;
$$

{\em (ii)} The automorphisms $(\f^* \rho)^2,$
where $\rho$ is any of $1$-permutations described
in {\em (i)}, are all $2$-transvections.
\end{Lem}

\begin{proof}
Let $\rho$ be a $1$-permutation satisfying the
conditions from (i). First note that due to
the proper commutativity with
$\theta^*,$  the restriction of $\rho$ to
$A^+_{\theta^*}$ must be equal to that one of $\pi^*.$

We denote by $R$ the matrix of the restriction
of $\rho$ on $B=A^-_{\theta^*}$
in the above described basis $\{e_1,e_2\}.$

Since the condition $(\pi^* \rho)^3=\id$
can be rewritten as
$$
\pi^* \rho \pi^* =\rho \pi^* \rho,
$$
we have
\begin{equation}
\begin{pmatrix}
0 & 1\\
1 & 0
\end{pmatrix}
R
\begin{pmatrix}
0 & 1 \\
1 & 0
\end{pmatrix}
=
R
\begin{pmatrix}
0 & 1 \\
1 & 0
\end{pmatrix}
R.
\end{equation}
Let now
$$
R =
\begin{pmatrix}
a & b\\
c & -a
\end{pmatrix},
$$
where $a,b,c \in \Z$ (the trace
of $R$ should be equal to zero
like the trace of any non-central
involution in $\gl{2,\Z},$ Theorem \ref{HRCanForm}).
It follows from (\theequation) that
\begin{align}
-a &= a(b+c),\\
c  &=b^2-a^2, \nonumber\\
b  &=c^2-a^2. \nonumber
\end{align}

According to (\theequation), there are two
cases for study: $a=0$ and $b+c+1=0.$

In the first case we have that $b=c=1$
and then $\rho=\pi^*,$ which is impossible.

The second case: we use the condition
$\det R =-1$ ($\rho$ is a conjugate
of $\pi^*$). Then
$$
\det R=-1 =-a^2-bc=-a^2-b(-b-1)
$$
or
$$
a^2=b^2+b+1.
$$
The only $b \in \Z$ for which
the number $b^2+b+1$ is
a square are $b=0,-1.$

Thus, there are indeed at most four possibilities for $R$:
\begin{equation}
R=
\begin{pmatrix}
e & 0 \\
-1    & -e
\end{pmatrix},
\begin{pmatrix}
e & -1 \\
0     & -e
\end{pmatrix},
\end{equation}
where $e=\pm 1.$ One easily verifies that for all four
$1$-permutations $\rho$ that correspond to the
matrices in (\theequation) and such that $\pi^* c
=\rho c$ for all $c \in A^+_{\theta^*},$ the conditions
from (i) of the Lemma are true.

The statement in (ii) is now a consequence of the
following observations:
\begin{align*}
\left[
\begin{pmatrix}
1 & 0 \\
0 & -1
\end{pmatrix}
\begin{pmatrix}
e & 0 \\
-1    & -e
\end{pmatrix}
\right]^2 =
\begin{pmatrix}
e^2 & 0 \\
2e  & e^2
\end{pmatrix}\\
\left[
\begin{pmatrix}
1 & 0 \\
0 & -1
\end{pmatrix}
\begin{pmatrix}
e & -1 \\
0     & -e
\end{pmatrix}
\right]^2
=
\begin{pmatrix}
e^2 & -2e\\
0   & e^2
\end{pmatrix},
\end{align*}
where $e=\pm 1.$
\end{proof}

\begin{Lem} \label{Def-of-2m-Trans}
The family of all $2m$-tranvections
{\em(}where $m$ runs over $\N${\em)} is $\varnothing$-definable
in $\aut A.$
\end{Lem}

\begin{proof}
We shall continue to use the parameters picked up above.  One more
parameter will be serviceable, however:  a
$2$-transvection $\tau^*,$ one of the four
$2$-transvections described in Lemma \ref{4-1-Perms} (ii).

Let us consider the
set $S$ of automorphisms $\{\f^* \rho\},$ where
$\rho$ is an extremal involution or
a $1$-permutation properly commuting
with $\theta^*.$ If the matrix
of the restriction of $\tau^*$ on $B$
is, for instance,
$$
\begin{pmatrix}
1 & 2 \\
0 & 1
\end{pmatrix}
$$
then only those elements from $S$
commute with $\tau^*$ whose restrictions
on $B$ have matrices
\begin{equation}
\begin{pmatrix}
e & b \\
0  & e
\end{pmatrix}
\end{equation}
where $e =\pm 1$ and $b \in \Z.$ The squares of the matrices of the form
(\theequation) are matrices
$$
\begin{pmatrix}
1 & 2b \\
0 & 1
\end{pmatrix}.
$$
Thus the set consisting of squares of elements of $S$
is a set that for each natural number $m$ contains a
$2m$-transvection. This implies that a suitably chosen existential
formula defines the $2m$-transvections.
\end{proof}

\begin{Lem} \label{MutSubgr}
Two distinct extremal involutions $\f_1,\f_2$ have
the mutual {\em(}eigen{\em)} subgroup, that is, either
$$
A^+_{\f_1}=A^+_{\f_2}, \text{ or }  A^-_{\f_1}=A^-_{\f_2}.
$$
if and only if the product $\f_2\f_1$
is a $2m$-tranvection for some non-zero
natural $m.$
\end{Lem}

\begin{proof}
Assume that subgroups $A^-_{\f_1}$ and
$A^-_{\f_2}$ are generated by
unimodular elements $x_1$ and $x_2$
respectively, and write $B_1$ and
$B_2$ for $A^+_{\f_1}$ and $A^+_{\f_2}.$
Let also $\tau$ denote
the product $\f_2\f_1.$

($\Leftarrow$). Suppose $B_1 \ne B_2.$ Then the
intersection $B_1 \cap B_2$ is of corank $2.$ The
fixed-point subgroup $C$ of $\tau$ has corank $1$ and
contains (a direct summand of $A$) $B_1 \cap B_2$;
then there is a unimodular element $y \in A$ such that
$$
C =\str y \oplus (B_1 \cap B_2).
$$
We have $\f_2 \f_1 y =y$ and then
$$
\f_1 y - y =\f_2 y -y.
$$
The above element is non-zero, since
otherwise $y \in B_1 \cap B_2.$ Thus
$\str{x_1} \cap \str{x_2} \ne 0,$ or
$\str{x_1}=\str{x_2},$ since both
$x_1,x_2$ are unimodular.

($\Rightarrow$). (i) Suppose that $B_1 \ne B_2,$ but $\str{x_1}
=\str{x_2}.$ Since $B_1 \cap B_2$ is a direct summand
of $A$ of corank $2,$ then for some unimodular $z$
$$
B_1 = \str z \oplus (B_1 \cap B_2);
$$
the element $z$ can be expressed
as $m x_2+b_2,$ where $m \in \Z$ and $b_2 \in B_2.$
We then have
$$
\tau z = \f_2 \f_1 z=\f_2 z=\f_2(m x_2+b_2) = -m x_2+b_2 =
m x_2+b_2 -2m x_2=z-2mx_2.
$$
Taking into account that $\tau x_2=x_2,$ we
see that $\tau$ is a $2m$-transvection.

(ii) Suppose that $B=B_1=B_2$ and
$\str{x_1} \ne \str{x_2}.$ The element
$x_1$ can be then written as
$$
x_1 = e x_2 + b=e x_2 +m c,
$$
where $b=m c$ is an element
of $B$ and $c$ is a unimodular. Hence
$$
\tau x_1 =\f_2\f_1 x_1 =\f_2(-e x_2-m c) =e x_2-m c=
x_1 -2m c
$$
and $\tau$ is a $2m$-transvection.
\end{proof}

\begin{Prop}
Let $\cD^1(A)$ be the family of all direct summands of
$A$ having rank or corank one. Then the action of the
group $\aut A$ on the family $\cD^1(A)$ is
first-interpretable in $\aut A$ without parameters.
\end{Prop}

\begin{proof}
In view of Lemma \ref{1-2-4-Def}, Lemma
\ref{Def-of-2m-Trans} and Lemma \ref{MutSubgr} all we
have to do is to explain when two pairs of extremal
involutions $(\f_1,\f_2)$ and $(\psi_1,\psi_2)$ both having
mutual subgroups determine the same direct summand of
$A$. It is easy: we just say that
for all $i,j$ either $\f_i=\psi_j,$
or $\f_i\psi_j$ is a $2m$-transvection.
\end{proof}

In  the conclusion of the section we present a purely algebraic
observation due to Oleg Belegradek who had found it while
reading the first draft of the paper.

\begin{Prop}
Let $A_1,A_2$ be infinitely generated
free abelian groups. The groups
$\aut{A_1}$ and $\aut{A_2}$ are
isomorphic if and only
if the cardinals $\rank A_1$ and $\rank A_2$
are equal.
\end{Prop}

\begin{proof}
Let $A$ be an infinitely generated free abelian group.
It is easy to show that the cardinality of any maximal
family of pairwise commuting 1-involutions in ${\rm
Aut}(A)$ is equal to rank of $A$. Since, by Lemma
1.4, the 1-involutions are $\varnothing$-definable in
${\rm Aut}(A)$ uniformly in $A$, and isomorphisms
preserve first-order formulae, the result follows.
\end{proof}

\section{Definability of the congruence subgroup of level two}

Let $m > 1$ be a natural number. Write $\Gamma_m(A)$
for the subgroup of $\aut A$ consisting of the
automorphisms of $A$ that act trivially (in the
natural way) on the group $A/mA.$ The subgroups
$\Gamma_m(A)$ are natural analogues of the principal
congruence subgroups of the groups $\spl{n,\Z}.$

We are going to prove $\varnothing$-definability of
the subgroup $\Gamma_2(A),$ the principal congruence
subgroup of $\aut A$ of level two. As it has
been said in the Introduction this will imply
a possibility of first-order interpretation
in $\aut A$ of the general linear group of
the vector space $A/2A$ over the field $\Z_2.$

\begin{Th} \label{Def-of-Gamma2A}
The subgroup $\Gamma_2(A)$
is $\varnothing$-definable in $\aut A.$
\end{Th}

\begin{proof} We shall use properties of the
group $\spl{3,\Z}$ and with this idea in mind we are
going to fix somehow some three direct
summands of rank one in $A.$ To achieve that we use
certain definable parameters. First, we take three
pairwise commuting extremal involutions
$\f_1^*,\f_2^*,\f^*_3$ in $\aut A$ such that any
product $\f_i^* \f_j^*$, where $i \ne j$ is a
$2$-involution. There exists a basis $\cB$ of $A$ in
which $\f_1^*,\f_2^*,\f_3^*$ are all diagonalizable.
Let $e_i$ denote the element of $\cB$ that $\f_i^*$
($i=1,2,3$) sends to the opposite.

Second, we need two $1$-permutations
$\pi_1^*$ and $\pi_2^*$ to provide
a suitable action on $\{e_1,e_2,e_3\}$;
our requirements on $\pi_1^*$ and
$\pi_2^*$ are therefore as follows:
\begin{itemize}
\item[(i)] $\pi_1^* \f_1^* \pi_1^* =\f_2^*$
and $\pi_1^*$ commutes with $\f_2^*$;

\item[(ii)] $\pi_2^* \f_1^* \pi_2^* =\f_3^*$
and $\pi_1^*$ commutes with $\f_3^*$;

\item[(iii)] $\pi_1^*$ and $\pi_2^*$
are conjugate and their product
is of order three.
\end{itemize}

In the following statement we simultaneously introduce
and characterize some transvections we are going to
deal with.

\begin{claim-num}
The {\em elementary transvections} which act trivially
on $\cB \setminus \{e_1,e_2,e_3\}$
and whose matrices in $\{e_1,e_2,e_3\}$
{\em(}more precisely, matrices of
the corresponding restrictions{\em)}
are of the form $E+n E_{ij},$ where
$1 \le i,j \le 3,$ $i \ne j$ and
$E_{ij}$ are the matrix units, are definable with parameters $\f_1^*,\f_2^*,\f_3^*$
and $\pi_1^*,\pi_2^*.$
\end{claim-num}

We choose a $2$-transvection
$\tau_1^*,$ one of the four $2$-transvections
that satisfy the condition (ii) of Lemma
\ref{4-1-Perms} for the $2$-involution
$\theta_1^*=\f_1^* \f_2^*$ and the $1$-permutation
$\pi_1^*.$ Without loss of generality we may suppose
that the matrix of $\tau_1^*$ in
$\{e_1,e_2,e_3\}$ is
$$
\begin{pmatrix}
1 & 2 & 0 \\
0 & 1 & 0 \\
0 & 0 & 1
\end{pmatrix}.
$$

It is easy to see that among the automorphisms $\f_1^*
\rho,$ where $\rho$ is either an extremal involution,
or a $1$-permutation properly commuting with
$\theta_1^*$ there are exactly four automorphisms
whose square is $\tau_1^*.$
The reason is that there are two solutions
to the matrix equation
$$
X^2 =
\begin{pmatrix}
1 & 2 \\
0 & 1
\end{pmatrix}
$$
in $\spl{2,\Z},$ namely,
$$
X =\pm
\begin{pmatrix}
1 & 1 \\
0 & 1
\end{pmatrix}
$$
and that any automorphism properly commuting
with $\theta_1^*$ must act on its fixed-point
subgroup, say, $C,$ either as the
$\id_C,$ or $-\id_C.$ Let us denote
the said four automorphisms by $\s_1,\s_2,\s_3,\s_4$ and let us further agree that $\s_1$ is
the {\it only} transvection among
the automorphisms $\s_i.$

The matrices of the automorphisms
$\s_i$ in the basis $\{e_1,e_2,e_3\}$
are
$$
\pm
\begin{pmatrix}
1 & 1 & 0 \\
0 & 1 & 0 \\
0 & 0 & 1
\end{pmatrix},
\pm
\begin{pmatrix}
-1 & -1 & 0 \\
0 & -1 & 0 \\
0 & 0 & 1
\end{pmatrix}.
$$
(the reader may as well imagine the diagonals
of the matrices stretched up to infinity filled
with units, but there is actually no need in that,
since already three coordinates do the
job.)

Let $\s$ be one of our automorphisms
$\s_i.$ We consider the conjugate $\s'=\pi \s \pi\inv$
of $\s$ by the automorphism $\pi=\pi_2^* \pi_1^*.$
Then the matrix of the commutator $[\sigma,\sigma']=
\s \s' \s\inv \s'{}\inv$
is either the matrix
$$
\begin{pmatrix}
1 & 0 & 1 \\
0 & 1 & 0 \\
0 & 0 & 1 \\
\end{pmatrix},
\text{ or }
\begin{pmatrix}
1 & 2 & -3 \\
0 & 1 & -2 \\
0 & 0 & 1
\end{pmatrix}.
$$

Thus only in the case when $\s=\s_1$ we have
the commutator $[\s,\s']$ {\it conjugate}
to $\s.$ Really, as to the automorphisms
$\s_2,\s_3,\s_4$ they all have eigen value
$-1,$ while none of the commutators $[\s_i,\s_i']$
with $i=2,3,4$ has this eigen value. Summing
up, we see that $\s_1,$ a $1$-transvection,
is definable over the chosen parameters.

Like in the proof of Lemma \ref{Def-of-2m-Trans}
we see that the elementary transvections whose
matrices in $\{e_1,e_2,e_3\}$ are
\begin{equation} \label{2m-transvs}
\begin{pmatrix}
1 & 2m & 0 \\
0 & 1 & 0 \\
0 & 0 & 1
\end{pmatrix}
\end{equation}
are definable with parameters $\f_1^*,\f_2^*,\pi_1^*$
and $\tau_1^*.$ Then elementary transvection
with the matrices
$$
\begin{pmatrix}
1 & n & 0 \\
0 & 1 & 0 \\
0 & 0 & 1
\end{pmatrix}
$$
are also definable with the parameters $\f_1^*,\f_2^*,\pi_1^* ,\pi_2^*,$
$\tau_1^*,$ since they are none the other than
either the transvections with matrices \eqref{2m-transvs},
or the products of the transvections with \eqref{2m-transvs}
and the elementary transvection with the matrix
$$
\begin{pmatrix}
1 & 1 & 0 \\
0 & 1 & 0 \\
0 & 0 & 1
\end{pmatrix}
$$
which is now known to be definable over
$\f_1^*,\f_2^*,\pi_1^*,\pi_2^*,\tau_1^*.$ The other
required elementary transvections are conjugates of
the tranvections with (\theequation) by suitable
automorphisms acting on $\{e_1,e_2,e_3\}$
as permutations, definable products of $\pi_1^*$ and
$\pi_2^*.$ Claim 1 is proved.

Let us note in passing that definability
of $1$-transvections with definable
parameters we have just proved
immediately implies the following
proposition.

\begin{Prop} Let $A$ be an infinitely
generated free abelian group. Then

{\em (i)} the family of all transvections
is $\varnothing$-definable in $\aut A;$

{\em (ii)} Let $m \ge 1$ be a natural number. The
family of all $m$-transvections is \linebreak
$\varnothing$-definable in $\aut A.$ \end{Prop}

Next is the construction of some set which is
contained in $\Gamma_2(A)$ and which is definable with
our parameters.

\begin{claim-num}
There is a set $D$ definable with parameters
$\f_1^*,$ $\f_2^*,$ $\pi_1^*,$ $\pi_2^*,$ $\tau_1^*,$ $\tau_2^*$
such that

{\em (i)} the automorphisms from $D$ act trivially
on $\cB\setminus \{e_1,e_2,e_3\}$ and their
matrices in $\{e_1,e_2,e_3\}$ are congruent
modulo $2$ to the identity matrix;

{\em (ii)} $D$ contains all automorphisms with {\em (i)}
whose matrices in $\{e_1,e_2,e_3\}$ are of the
form
$$
\begin{pmatrix}
a & b & 0 \\
c & d & 0 \\
0 & 0 & 1
\end{pmatrix}
$$
where
$$
a \equiv d \equiv 1\,(\operatorname{mod} 2) \text{ and }
b \equiv c \equiv 0\,(\operatorname{mod} 2).
$$
\end{claim-num}

The argument is based upon the remarkable observation
made in the paper \cite{CarKel} by Carter and Keller:
\begin{quote}
each matrix of the form
$$
\begin{pmatrix}
a & b & 0 \\
c & d & 0 \\
0 & 0 & 1
\end{pmatrix}
$$
from the (matrix) group $\spl{3,\Z}$ is
a product of at most 41 elementary transvections.
\end{quote}

Suppose that $t_1,\ldots, t_{41}$ are elementary
transvections, matrices from $\spl{3,\Z}.$ One
corresponds to the product
$$
t_1 t_2 \ldots t_{41}
$$
a sequence
\begin{equation}
(\av t_1,\av t_2,\ldots, \av t_{41})
\end{equation}
where $\av t$ is the image of $t$ in $\spl{3,\Z_2}$
under the natural homomorphism $\spl{3,\Z} \to \spl{3,\Z_2}.$
There are of course finitely many sequences
of the form (\theequation). Some of them determine
the identity matrix in $\spl{3,\Z_2},$ some do not;
we appreciate the former sequences, say `good'
ones. Clearly, the image $\av t$ of an {\it elementary}
transvection $t$ is trivial in $\spl{3,\Z_2}$
if and only if $t$ is a square of an elementary transvection in $\spl{3,\Z}.$
So the fact that a sequence $(\av t_1,\av t_2,\ldots, \av t_{41})$
is `good' can be translated into a disjunction
of statements each of which says for every $i=1,\ldots,41$
that the $i$th transvection $t_i$
is or is not a square.

Having the elementary transvections with
respect to the basis $\{e_1,e_2,e_3\}$
(this time automorphisms of $A$) definable
in $\aut A$ with the parameters introduced above,
we may realize the above considerations
for the group $\aut A.$ This completes
the proof of Claim 2.

Let now $\chi(\av v)$ be a first-order formula that
describes the parameters $\f_1^*,$ $\f_2^*,$
$\pi_1^*,$ $\pi_2^*,$ $\tau_1^*,$ $\tau_2^*.$ Suppose
that $\av \f$ is any tuple of elements of $\aut A$
that satisfies $\chi$; we then denote by $D(\av \f)$
the family of automorphisms constructed over $\av\f$
in the same way as $D$ is constructed over our
parameters.

\begin{claim-num}
The following are equivalent:

{\em (a)} $\sigma \in \aut A$ is an element
of $\Gamma_2(A);$

{\em (b)} there is a direct summand $B$ of $A$
of rank or corank $1$ such that for every
direct summand $C$ isomorphic via some
automorphism from $\aut A$ to $B$
there exist a tuple $\av \f$ satisfying
$\chi$ and $\rho \in D(\av \f)$ with
$$
\sigma C =\rho C.
$$
\end{claim-num}

Let consider the implication (b) $\Rightarrow$ (a).
Suppose that the direct summand $B$ mentioned
in (b) is of rank one and $e$ a unimodular
element of $A.$ Then for suitable
parameters $\av\f$ there is $\rho \in D(\av \f)$
$$
\s\str e =\rho \str e.
$$
By Claim 2 the set $D(\av \f)$ is contained
in $\Gamma_2(A)$ and hence
$$
\s e =\pm \rho e \equiv \pm e \equiv e (\operatorname{mod} 2A).
$$
It then follows that $\s \in \Gamma_2(A).$

Suppose now that $B$ is of corank 1.
Let $e$ be a unimodular element of $A$
and let $\{e,e_0,e_1,\ldots,e_n,\ldots\}$
be a basis of $A.$ According to the condition $\sigma$
moves the direct summand
$$
C_0 = \str{e,e_1,e_2,\ldots,e_n,\ldots}
$$
exactly as some $\rho \in \Gamma_2(A)$ does:
$$
\sigma C_0 =\rho C_0.
$$
This implies that $\sigma e$ is
congruent modulo $2A$ to
some element of $C_0$:
\begin{equation}
\sigma e \equiv ke+ k_1 e_1+k_2 e_2+\ldots+k_n e_n +\ldots (\operatorname{mod} 2A).
\end{equation}
The same argument can be applied to the subgroup
$$
C_1=\str{e,e_0,e_2,\ldots,e_n,\ldots}
$$
of which $e$ is also a member; this leads to
$$
\sigma e \equiv le+l_0 e_0 +l_2 e_2+\ldots+l_n e_n+\ldots (\operatorname{mod} 2A).
$$
One deduces then that
$$
(k-l)e-l_0 e_0+k_1 e_1+(k_2-l_2) e_2+\ldots+(k_n-l_n) e_n+\ldots \equiv 0 (\operatorname{mod} 2A).
$$
The images of $e,e_0,e_1,e_2,\ldots$ under the natural
homomorphism $A \to A/2A$ must be linearly independent
over $\Z_2$ and therefore
$$
l_0 \equiv k_1 \equiv 0\, (\operatorname{mod} 2).
$$
Continuing in a similar fashion, we see that
all (non-zero) coefficients $k_i$ in (\theequation)
are even; the coefficient $k$ must therefore be odd.
Thus $\s$ is in $\Gamma_2(A),$ as required.

The implication (a) $\Rightarrow$ (b).
Suppose that $\sigma \in \Gamma_2(A)$
and $e$ is a unimodular element of $A.$
Then for a basis $\{e,e_0,e_1,\ldots,e_n,\ldots\}$
of which $e$ forms a part we have
$$
\sigma e =e +2(ke +\sum_i k_i e_i).
$$
Suppose that $s$ is the greatest
common divisor of non-zero
elements $k_i.$ Then
$$
\sigma e = (1+2k) e+ 2s(\sum_i k_i' e_i).
$$
Clearly, $\gcd(1+2k,2s)=1$ (since
$\sigma e$ is unimodular) and
the element $g=\sum_i k_i' e_i$
is unimodular. If so, there
are $b,d \in \Z$ such that the
matrix
$$
\begin{pmatrix}
1+2k & b & 0 \\
2s   & d & 0 \\
0    & 0 & 1
\end{pmatrix}
$$
from $\spl{3,\Z}$ is congruent to the identity
matrix modulo $2.$ This implies
that there exist a tuple $\av \f$ satisfying
$\chi$ and some $\rho \in D(\av \f)$ such that
$$
\sigma \str e =\rho \str e.
$$
Claim 3 is proved.

Since we know how to interpret in $\aut A$ by means of
first-order logic the direct summands of $A$ of
rank/corank $1,$ the conditions in (ii) of Claim 3
are easily translated into first-order formulae. The
proof of Theorem \ref{Def-of-Gamma2A} is now
completed.
\end{proof}

\begin{Rem} {\em Very recently Bardakov
proved that the principal congruence
subgroups of the groups $\spl{n,\Z},$ where $n \ge 3$
all have finite width with respect to elementary
transvections (unpublished; personal communication).
Recall that the {\it width} of a group $G$ relative to a generating
set $S$ with $S\inv=S$ is either the minimal
natural number $k$ such that every
element of $G$ is a product of at most
$k$ elements of $S,$ or $\infty$ otherwise.

The result by Bardakov could be used then
to simplify the proof of Theorem \ref{Def-of-Gamma2A}.}
\end{Rem}

\begin{Th} \label{Int-Set-Theo}
Let $A$ be an infinitely generated
free abelian group. Then the group
$\aut A$ first-order interprets
the second-order theory $\Theo_2(|A|),$
uniformly in $A.$
\end{Th}

\begin{proof}
The proof is based on Theorem \ref{Def-of-Gamma2A}
and the following important theorem from the paper \cite{BrMa}
by Bryant and Macedonska.

\begin{theo}
Let $F$ be a free group of infinite rank
and let $V$ be a characteristic subgroup
of $F$ such that $F/V$ is nilpotent.
Then every automorphism of $F/V$
is induced by an automorphism of
$F.$
\end{theo}

Let $A$ stand for the free abelian group $F/[F,F].$ As
a corollary of the result by Bryant--Macedonska we have that the natural homomorphism
$$
\mu : \aut A \to \aut{A/2A}
$$
(induced by the natural homomorphism $A \to A/2A$) is surjective. Indeed, according to
the Theorem, the natural homomorphisms
$$
\mu_1 : \aut F \to \aut A \text{ and }
\mu_2 : \aut F \to \aut{A/2A}
$$
are both surjective. On the other
hand,
$$
\mu_2 = \mu \circ \mu_1,
$$
and then $\mu$ must be surjective, too.

Adding this to the fact that
$\Gamma_2(A),$ the kernel of
$\mu,$ is $\varnothing$-definable
in $\aut A,$ we get that
the group $\aut A$ first-order
interprets the group $\aut{A/2A}$:
$$
\aut A/\ker \mu =\aut A/\Gamma_2(A) \cong \aut{A/2A}.
$$

The group $\aut{A/2A}$ is the general
linear group of the vector space
$A/2A$ over the field $\Z_2.$
On the other hand, the general linear group $\gl V$ of
a infinite-dimensional vector space $V$ over a field
$D$
first-order interprets $\Theo_2(\dim_D V),$ see
\cite[Theorem 11.4]{ToAPAL}. Therefore the elementary
theory of the group $\aut{A/2A}$
first-order interprets
the second-order theory
$$
\Theo_2(\dim_{\Z_2} A/2A)=\Theo_2(|A|),
$$
and the result follows.
\end{proof}

\begin{cor}
Let $A_1,A_2$ be infinitely generated
free abelian groups. The groups
$\aut{A_1}$ and $\aut{A_2}$ are
elementarily equivalent if and only
if the cardinals $|A_1|$ and $|A_2|$ {\em(}viewed
as sets with no structure{\em)}
are second-order equivalent.
\end{cor}

\begin{proof}
The necessity part is a consequence of Theorem \ref{Int-Set-Theo}.
To prove the converse, one syntactically interprets
in the second-order theory $\Theo_2(\vk),$
where $\vk$ is an infinite cardinal,
the elementary theory of the automorphism
group of a free abelian
group with $\vk$ as the domain (rather easy;
cf. \cite[Theorem 4.1]{ToJLM2} where a similar
interpretation is done in quite full detail for the case of
the elementary theory of the automorphism
group a free group
over $\vk.$)
\end{proof}


\begin{thebibliography}{99}


\bibitem{BrMa} R.~Bryant, O.~Macedonska, {\it
Automorphisms of relatively free nilpotent groups of
infinite rank}, J. Algebra {\bf 121} (1989), 388--398.


\bibitem{CarKel} D.~Carter, G.~Keller, {\it Elementary
expressions for unimodular matrices}, {Commun.
Algebra}, {\bf 12} (1984), 379--389.

\bibitem{HuaRei} L.~K.~Hua, I.~Reiner. {\it Automorphisms of
the unimodular group}, {Trans. Amer. Math.
Soc.} {\bf 71} (1951), 331--348.

\bibitem{O'M} O.~T.~O'Meara. {\it Lectures on linear
 groups}, Amer. Math. Soc., Providence, RI, 1974.


\bibitem{Sh1} S.~Shelah, {\it First-order theory of permutation
groups}, {Israel. J. Math.} {\bf 14} (1973), 149--162.


\bibitem{ShCler} S.~Shelah, {\it Interpreting set
theory in the endomorphism semi-group of a free
algebra or in a category}, {Annales Scientifiques de L'universite
de Clermont} {\bf 13} (1976), 1--29.

\bibitem{ShMatJap} S.~Shelah, {\it On what I do not
understand {\em(}and have something to say{\em)},
model theory}, Math. Japon. {\bf 51} (2000), 329--377.


\bibitem{ToAPAL} V.~Tolstykh, {\it Elementary equivalence of
infinite-dimensional classical groups}, {Ann. Pure Appl.
Logic} {\bf 105} (2000), 103--156.


\bibitem{ToJLM2} V.~Tolstykh, {\it Set theory is interpretable in the
automorphism group of an infinitely generated free
group}, {J. London Math. Soc.} {\bf 62} (2000),
16--26.

\bibitem{ToContMat} V.~Tolstykh, {\it On the logical strength of
the automorphism groups of free nilpotent groups},
Contemp. Math. vol. 302, AMS, Providence, 2002, 113--120.



\bibitem{ToCamb} V.~Tolstykh, {\it Free two-step nilpotent groups
whose automorphism group is complete}, Math. Proc.
Cambridge Philos. Soc. {\bf 131} (2001), 73--90.




\end{thebibliography}
\end{document}